\documentclass[12pt,a4paper]{article}
\usepackage{amsfonts}

\newtheorem{theorem}{Theorem}[section]

\newtheorem{prop}{Proposition}[section]
\newtheorem{cor}{Corollary}[section]

\newcommand{\LP}{L\'{e}vy process}
\newcommand{\R}{\mathbb{R}}
\newcommand{\Rd}{\mathbb{R}^{d}}
\newcommand{\nN}{n \in \mathbb{N}}
\newcommand{\N}{\mathbb{N}}
\newcommand{\C}{\mathbb{C}}
\newcommand{\E}{\mathbb{E}}

\newcommand{\cadlag}{c\`adl\`ag}
\newcommand{\ds}{\displaystyle}
\newcommand{\g}{\textbf{g}}

\newcommand{\bean}{\begin{eqnarray*}}
\newcommand{\eean}{\end{eqnarray*}}

\newcommand{\la}{\langle}
\newcommand{\ra}{\rangle}

\newcommand{\Z}{\mathbb{Z}_{+}}

\newcommand{\at}{\texttt{a}^{*}}

\date{}

\begin{document}

\title{Aspects of Recurrence and Transience for L\'{e}vy Processes in
Transformation Groups and Non-Compact Riemannian Symmetric Pairs}

\author{ David Applebaum, \\ School of Mathematics and Statistics,\\ University of
Sheffield,\\ Hicks Building, Hounsfield Road,\\ Sheffield,
England, S3 7RH\\ ~~~~~~~\\e-mail: D.Applebaum@sheffield.ac.uk }

\maketitle

\begin{abstract}\noindent
We study recurrence and transience for L\'{e}vy processes induced
by topological transformation groups. In particular the
transience-recurrence dichotomy in terms of potential measures is
established and transience is shown to be equivalent to the
potential measure having finite mass on compact sets when the
group acts transitively.  It is known that all bi-invariant
L\'{e}vy processes acting in irreducible Riemannian symmetric
pairs of non-compact type are transient. We show that we also have
``harmonic transience'', i.e. local integrability of the inverse
of the real part of the characteristic exponent which is
associated to the process by means of Gangolli's
L\'{e}vy-Khinchine formula.

\vspace{5pt}

\noindent  {\it Key words and phrases.}~Transformation group, Lie
groups, L\'{e}vy process, recurrence, transience, potential
measure, convolution semigroup, spherical function, Gangolli's
L\'{e}vy-Khintchine formula, transient Dirichlet space, harmonic
transience.

\vspace{5pt}

\noindent {\it AMS 2000 subject classification.}~60B15, 60G51,
31C25, 60J45, 37B20.

\end{abstract}

\section{Introduction}

L\'{e}vy processes in Lie groups have recently attracted
considerable interest and a recent monograph \cite{Liao} is
dedicated to their investigation. The purpose of this paper is to
develop some theoretical insight into the recurrence and
transience of such processes. L\'{e}vy processes are often
considered as natural continuous-time generalisations of random
walks. It is well known that recurrence and transience of random
walks in groups is intimately related to volume growth in the
group (see \cite{Bal}, \cite{BLP}, \cite{GKR}). From another point
of view a path continuous L\'{e}vy process on a group is a
Brownian motion with drift and there has been a great deal of work
on transience/recurrence of Brownian motion in the more general
context of processes on Riemannian manifolds. For a nice survey
see \cite{Gr}.

The classic analysis of transience and recurrence for L\'{e}vy
processes on locally compact abelian groups involves a subtle
blend of harmonic analytic and probabilistic techniques. This work
was carried out by Port and Stone \cite{PS} who showed that a
necessary and sufficient condition for transience is local
integrability of the real part of the inverse of the
characteristic exponent with respect to the Plancherel measure on
the dual group. A non-probabilistic version of this proof was
given by M.It\^{o} \cite{It} - see also \cite{HT} and section 6.2
in \cite{Hey}. It relies heavily on reduction to the case where
the group is of the form $\R^{d} \times \mathbb{Z}^{n}$.


Harmonic analytic methods may, at least in principle, be applied
to the non-commutative case if we work with Gelfand pairs $(G,K)$
(or more generally, with hypergroups for which we refer readers to
section 6.3 of \cite{BH}.) Here we can take advantage of the
existence of spherical functions to develop harmonic analysis of
probability measures in the spirit of the abelian case - indeed in
the important case where $G$ is a connected Lie group and $K$ is a
compact subgroup (so that $G/K$ is a symmetric space), a L\'{e}vy
-Khintchine type formula which classifies bi-invariant L\'{e}vy
processes in terms of their characteristic exponent was developed
by Gangolli \cite{Gang}. For further developments of these ideas
see \cite{App0}, \cite{LW} and the survey article \cite{Hey1}. An
important approach to establishing the transience of a Markov
process is to prove that the associated Dirichlet form gives rise
to a transient Dirichlet space. A comprehensive account of this
approach can be found in section 1.5 of \cite{FOT}. Bi-invariant
Dirichlet forms associated to Gelfand pairs $(G,K)$ were first
studied in a beautiful paper by C.Berg \cite{Berg}. He was able to
establish that if the Dirichlet space is transient then the
inverse of the characteristic exponent is locally integrable with
respect to Plancherel measure on the space of positive definite
spherical functions. However he was only able to establish the
converse to this result in the case where the group was compact or
of rank one. By using different techniques,  Berg and Faraut
\cite{BeFa} (see also the survey \cite{Hey2}) were able to show
that all bi-invariant L\'{e}vy processes associated to non-compact
irreducible Riemannian symmetric pairs are transient as are the
associated Dirichlet spaces in the symmetric case. One of the
goals of the current paper is to show that all such processes are
{\it harmonically transient} in the sense that the inverse of real
part of the characteristic exponent is locally integrable with
respect to Plancherel measure. Note that this is a slightly
stronger result than is obtained in the Euclidean case.

The organisation of this paper is as follows. First we study
L\'{e}vy processes in a quite abstract context, namely topological
transformation groups acting on a complete metric space. This
allows us to work with processes on the group and then study the
induced action on the space of interest. There are three key
results here. First we establish a recurrence/transience dichotomy
for the induced process in terms of potential measures of open
balls. This part of the paper closely follows the development
given in section 7.35 of Sato \cite{Sato} for Euclidean spaces.
Secondly we show that such processes are always recurrent when the
space is compact and thirdly we establish that transience of the
process is (given that the group acts transitively) equivalent to
the finiteness of the potential measure on compact sets. This last
result (well known in the abelian case) is important for bridging
the gaps between probabilistic and analytic approaches to
transience.

In the second half of the paper we specialise by taking the group
$G$ to be a non-compact semisimple group having a finite centre.
We fix a compact subgroup $K$  so that we can consider the action
on the symmetric space $G/K$. We consider symmetric
$K$-bi-invariant L\'{e}vy processes within this context. Using the
spherical transform, we establish pseudo-differential operator
representations of the Markov semigroup, its generator and the
associated Dirichlet form which may be of interest in their own
right (c.f. \cite{App2}.) We then establish the result on harmonic
transience as described above.

\bigskip

{\bf Notation} If $M$ is a topological space, ${\cal B}(M)$ is the
$\sigma$-algebra of Borel measurable subsets of $M$, $B_{b}(M)$ is
the Banach space of all bounded Borel measurable function on $M$
(when equipped with the supremum norm $|| \cdot ||$), $C_{0}(M)$
is the closed subspace of continuous functions on $M$ which vanish
at infinity and $C_{c}(M)$ is the space of continuous functions on
$M$ which have compact support. If ${\cal F}(M)$ is any real
linear space of real-valued functions on $M$, then ${\cal
F}_{+}(M)$ always denotes the cone therein of non-negative
elements. Throughout this article, $G$ is a topological group with
neutral element $e$. For each $\sigma \in G, l_{\sigma}$ denotes
left translation on $G$ and $l_{\sigma}^{*}$ is its differential.
If $m$ is a Haar measure on a locally compact group $G$, we write
$m(d\sigma)$ simply as $d\sigma$. When $G$ is compact, we always
take $m$ to be normalised. The {\it reversed measure}
$\widetilde{\mu}$ that is associated to each Borel measure $\mu$
on a topological group $G$ is defined by $\widetilde{\mu}(A): =
\mu(A^{-1})$ for each $A \in {\cal B}(G)$. If $f \in L^{1}(G,\mu)$
we will sometimes write $\mu(f) := \int_{G}fd\mu$. If $X$ and $Y$
are $G$-valued random variables defined on some probability
spaces, with laws $P_{X}$ and $P_{Y}$, respectively, then $X
\stackrel{d}{=} Y$ means $P_{X} = P_{Y}$. $\R^{+}:=[0, \infty)$.
We will use Einstein summation convention throughout this paper.

\section{Transformation Groups and L\'{e}vy Processes}

\subsection{Probability on Transformation Groups}

Let $G$ be a topological group with neutral element $e$, $M$ be a
topological space and $\Phi: G \times M \rightarrow M$ be
continuous. We say that $(G,M, \Phi)$ is a {\it transformation
group} if for all $m \in M$,
\begin{enumerate}
\item[(T1)] $\Phi(e,m) = m$. \item[(T2)] $\Phi(\sigma, \Phi(\tau,
m)) = \Phi(\sigma \tau, m)$, for all $\sigma, \tau \in G$.
\end{enumerate}

For fixed $m \in M$, we will often write $\Phi_{m}$ to denote the
continuous map from $G$ to $M$ defined by $\Phi_{m}: = \Phi(\cdot,
m)$.

The transformation group $(G,M, \Phi)$ is said to be {\it
transitive} if for all $m,p \in M$ there exists $\sigma \in G$
such that $p = \Phi(\sigma, m)$. In this case each mapping
$\Phi_{m}$ is surjective. A rich class of transitive
transformation groups is obtained by choosing a closed subgroup
$K$ of $G$ and taking $M$ to be the homogeneous space $G/K$ of
left cosets. In this case, we will write
     $$ \Phi_{\sigma}(p) := \Phi(\sigma, p) = \sigma \tau K,$$
where $p = \tau K$ for some $\tau \in G$. The canonical surjection
from $G$ to $G/K$ will be denoted by $\pi$.

Now let $(\Omega, {\cal F}, P)$ be a probability space and $X$ be
a $G$-valued random variable with law $p_{X}$. For each $m \in M$
we obtain an $M$-valued random variable $X^{\Phi,m}$ by the
prescription $X^{\Phi,m}: = \Phi_{m}(X)$. The law of $X^{\Phi,m}$
is $p_{X}^{\Phi,m}: = p_{X} \circ \Phi_{m}^{-1}$ and it is clear
that if $Y$ is another $G$-valued random variable then
$X^{\Phi,m}$ and $Y^{\Phi,m}$ are independent if $X$ and $Y$ are.
Sometimes we will work with a fixed $m \in M$ and in this case we
will write $X^{\Phi} = X^{\Phi,m}$. In the case where $M = G/K$ we
will always take $m = eK$ and write $X^{\pi}:= X^{\Phi,m}$.

\subsection{Metrically Invariant Transformation Groups}

For much of the work that we will carry out in this paper we will
need additional structure on the transformation group
$(G,X,\Phi)$. Specifically we will require the space $X$ to be
metrisable by a complete invariant metric $d$. In this case we
will say that $(G,X,\Phi)$ is {\it metrically invariant}. We
emphasise that invariance in this context is the requirement that
  $$ d(\Phi(\sigma,  m_{1}), \Phi(\sigma, m_{2})) = d(m_{1},
  m_{2}),$$ for all $m_{1}, m_{2} \in X, \sigma \in G$.

Two examples will be of particular relevance for our work:

\vspace{5pt}

{\bf Example 1}. {\it Compact Lie Groups}

\vspace{5pt}

Here $M$ = $G$ is a compact Lie group and $\Phi$ is left
translation. If $G$ is equipped with an Ad$(G)$-invariant positive
definite Riemannian metric, the associated distance function
inherits left-invariance. Since $G$ is compact, it is geodesically
complete and so is a complete metric space by the Hopf-Rinow
theorem (see e.g. \cite{Chav}, pp.\ 26--7).

Such a metric always exists when $G$ is also connected and
semi-simple. In this case we can, for example, take $B$ to be the
associated Killing form, then $-B$ is positive definite. We thus
obtain a left-invariant Riemannian metric $\langle \cdot, \cdot
\rangle$ via the prescription
$$
    \langle X_{\sigma}, Y_{\sigma} \rangle = -
    B(l_{\sigma^{-1}}^{*}(X_{\sigma}), l_{\sigma^{-1}}^{*}(Y_{\sigma})),
$$
for each $X_{\sigma}, Y_{\sigma} \in T_{\sigma}(G)$.

\vspace{5pt}

{\bf Example 2}. {\it Symmetric Spaces} \cite{Hel1}, \cite{Wo}

\vspace{5pt}

Let $G$ be a Lie group that is equipped with an involutive
automorphism $\Theta$ and let $K$ be a compact subgroup of $G$
such that $\Theta(K) \subseteq K$. $M = G/K$ is then a
$C^{\infty}$-manifold and we take $m = \pi(e) = eK$. Let $\g$ be
the Lie algebra of $G$. We may write $\g = \textbf{l} \oplus
\textbf{p}$ where $\textbf{l}$ and $\textbf{p}$ are the $+1$ and
$-1$ eigenspaces of $\Theta^{*}$, respectively and we have
$\pi^{*}(\textbf{p}) = T_{m}(M)$. Every positive definite
Ad$(K)$-invariant inner product defines a $G$-invariant Riemannian
metric on $M$ under which $M$ becomes a globally Riemannian
symmetric space (with geodesic symmetries induced by $\Theta$).
Integral curves of vector fields on $G$ project to geodesics on
$M$ and completeness follows by the same argument as in Example 1.

In particular if $G$ is semi-simple then $\textbf{l}$ =
$\textbf{k}$  where $\textbf{k}$ is the Lie algebra of $K$ and $\g
= \textbf{k} \oplus \textbf{p}$ is a Cartan decomposition of
$\textbf{g}$. Using the fact that any tangent vector $Y_{q} \in
T_{q}(M)$ (where $q = \tau K$ for $\tau \in G$) is of the form
$l_{\tau}^{*} \circ \pi^{*}(Y)$ for some $Y \in \textbf{p}$, the
required metric is given by

$$ \langle Y_{q}, Z_{q} \rangle =
    B(Y,Z),$$ where $Z_{q} = l_{\tau}^{*} \circ \pi^{*}(Z)$ and $B$ is the Killing form on
    $\textbf{g}$.

\vspace{5pt}

A third important class of examples is obtained by taking $M$ to
be an arbitrary geodesically complete Riemannian manifold and $G$
to be the group of all isometries of $M$.

\section{L\'{e}vy Processes in Groups}

Let $Z = (Z(t), t \geq 0)$ be a stochastic process defined on
$(\Omega, {\cal F}, P)$, and taking values in the topological
group $G$. The left increment of $Z$ between $s$ and $t$ where $s
\leq t$ is the random variable $Z(s)^{-1}Z(t)$.

We say that $Z$ is a \emph{L\'{e}vy process in $G$} if it
satisfies the following:
\begin{enumerate}
\item
    $Z$ has stationary and independent left increments,
\item
    $Z(0) = e$ (a.s.),
\item
    $Z$ is stochastically continuous,
i.e. $\lim_{s \downarrow t} P(Z(s)^{-1}Z(t) \in A) = 0 $ for all
$A \in {\cal B}(G)$ with $e \notin \bar{A}$ and all $t \geq 0$.
\end{enumerate}

Now let $(\mu_{t}, t \geq 0)$ be the law of the L\'{e}vy process
$Z$, then it follows from the definition that $(\mu_{t}, t \geq
0)$ is a weakly continuous convolution semigroup of probability
measures on $G$, where the convolution operation is defined for
probability measures $\mu$ and $\nu$ on $G$  to be the unique
probability measure $\mu * \nu$ such that
$$   \int_{G}f(\sigma)(\mu * \nu)(d\sigma) = \int_{G}\int_{G}f(\sigma \tau) \mu(d\sigma)\nu(d\tau), $$

for each $f \in B_{b}(G)$. So that in particular we have, for all
$s,t \geq 0$
\begin{equation}  \label{conv}
   \mu_{s+t} = \mu_{s} * \mu_{t}
   \qquad\textrm{and}\qquad
   \lim_{t \rightarrow 0} \mu_{t} = \mu_{0} = \delta_{e},
\end{equation}

where $\delta_{e}$ is Dirac measure concentrated at $e$, and the
limit is taken in the weak topology of measures. We obtain a
contraction semigroup of linear operators $(T_{t}, t \geq 0)$ on
$B_{b}(G)$ by the prescription

\begin{equation} \label{sem}
    (T_{t}f)(\sigma)  =  \E(f(\sigma Z(t)))
                    =  \int_{G}f(\sigma \tau)\mu_{t}(d\tau),
\end{equation}
for each $t \geq 0, f \in B_{b}(G), \sigma \in G$. The semigroup
is left-invariant in that $L_{\sigma}T_{t} = T_{t} L_{\sigma}$ for
each $t \geq 0, \sigma \in G$ where $L_{\sigma}f(\tau) =
f(\sigma^{-1}\tau)$ for all $f \in B_{b}(G), \tau \in G$.

Conversely, given any weakly continuous convolution semigroup of
probability measures $(\mu_{t}, t \geq 0)$ on $G$ on we can always
construct a L\'{e}vy process $Z = (Z(t), t \geq 0)$ such that each
$Z(t)$ has law $\mu_{t}$ by taking $\Omega$ to be the space of all
mappings from $\R^{+}$ to $G$ and ${\cal F}$ to be the
$\sigma$-algebra generated by cylinder sets. The existence of $P$
then follows by Kolmogorov's construction and the time-ordered
finite-dimensional distributions have the form \bean & &
P(Z(t_{1}) \in A_{1}, Z(t_{2}) \in A_{2}, \ldots , Z(t_{n}) \in
A_{n})\\ & = &  \int_{G} \int_{G} \cdots
\int_{G}1_{A_{1}}(\sigma_{1})1_{A_{2}}(\sigma_{1}\sigma_{2})
\cdots 1_{A_{n}}(\sigma_{1}\sigma_{2}\cdots
\sigma_{n})\mu_{t_{1}}(d\sigma_{1})\mu_{t_{2}-t_{1}}(d\sigma_{2})
\ldots \mu_{t_{n}-t_{n-1}}(d\sigma_{n}), \eean for all $A_{1},
A_{2}, \ldots A_{n} \in {\cal B}(G), 0 \leq t_{1} \leq t_{2} \leq
\ldots \leq t_{n} < \infty$. For a proof in the case $G = \Rd$,
see Theorem 10.4 on pages 55-7 of \cite{Sato}. The extension to
arbitrary $G$ is straightforward. The L\'{e}vy process $Z$ that is
constructed by these means is called a canonical \LP.

If $G$ is a locally compact group then $(T_{t}, t \geq 0)$ is a
left-invariant Feller semigroup in that
$$
    T_{t}(C_{0}(G)) \subseteq C_{0}(G)
    \qquad\textrm{and}\qquad
    \lim_{t \downarrow 0}||T_{t}f - f|| = 0
$$
for each $f \in C_{0}(G)$. The infinitesimal generator of $(T_{t},
t \geq 0)$ is denoted by ${\cal A}$. A characterisation of ${\cal
A}$ can be found in \cite{Born} ( see \cite{Hunt}, \cite{Liao} for
the Lie group case.) It follows from the argument on page 10 of
\cite{Liao} that if $(Z(t), t \geq 0)$ is a $G$-valued Markov
process with left-invariant Feller transition semigroup then it is
a \LP. Moreover if $G$ is separable and metrizable as well as
being locally compact then by Theorem 2.7 in Chapter 4 of
\cite{EK} p.169, the process has a \cadlag~modification (i.e. one
that is almost surely right continuous with left limits.)

Let $(G,M,\Phi)$ be a transformation group and $(Z(t), t \geq 0)$
be a L\'{e}vy process on $G$. Then for each $m \in M, t \geq 0$ we
define $Z^{\Phi,m}(t): = \Phi_{m}(Z(t))$. The $M$-valued process
$Z^{\Phi,m} = (Z^{\Phi,m}(t), t \geq 0)$ will be called a {\it
L\'{e}vy process on $M$ starting at $m$}. We obtain a contraction
semigroup of linear operators $(S_{t}^{\Phi}, t \geq 0)$ on
$B_{b}(M)$ by the prescription
$$  S_{t}^{\Phi}f(m) = T_{t}(f \circ \Phi_{m})(e),$$ for all $f
\in B_{b}(M), t \geq 0$. The case where $M$ is a homogeneous space
is discussed in \cite{Liao}, section 2.2. Suppose that we are
given a filtration $({\cal F}_{t}, t \geq 0)$ of ${\cal F}$ so
that $(Z(t), t \geq 0)$ is ${\cal F}_{t}$-adapted and that each
increment $Z(s)^{-1}Z(t)$ is independent of ${\cal F}_{s}$ for $0
< s < t < \infty$. In this case $Z$ is a Markov process (with
respect to the given filtration) having transition semigroup
$(T_{t}, t \geq 0)$. Moreover if $(G,M, \Phi)$ is transitive then
$Z^{\Phi,Z(0)}$ is a Markov process having transition semigroup
$(S_{t}^{\Phi}, t \geq 0)$. If $G$ is locally compact then
$(S_{t}^{\Phi}, t \geq 0)$ is a Feller semigroup whose generator
will be denoted ${\cal A}^{\Phi}$. We have
$$ \mbox{Dom}({\cal A}^{\Phi}) = \{f \in C_{0}(M); f \circ \Phi_{m}
\in \mbox{Dom}({\cal A})~\mbox{for all}~m \in M\}$$ and
$${\cal A}^{\Phi}f(m) = {\cal A}(f \circ \Phi_{m})(e),$$
for all $m \in M, f \in \mbox{Dom}({\cal A}^{\Phi})$.

\section{Criteria for Recurrence and Transience}

Throughout this section $(G,M,\Phi)$ will be a metrically
invariant transformation group and $d$ will denote the complete
metric on $M$, $Z$ will be a \LP~on $G$ and $Z^{\Phi,m}$ will the
associated L\'{e}vy process on $M$ starting at $m$.

We say that $Z$ is \emph{recurrent at $m$} if

$$
    \liminf_{t \rightarrow \infty} d(Z^{\Phi,m}(t),m) = 0\quad\textrm{(a.s.)},
$$
and \emph{transient at $m$} if
$$
    \lim_{t \rightarrow \infty} d(Z^{\Phi,m}(t), m) = \infty\quad\textrm{(a.s.)}.
$$

If $(G,M,\Phi)$ is also transitive, it follows easily from the
invariance of $d$ that if $Z$ is recurrent (respectively,
transient) at any given point of $M$ then it is recurrent
(respectively, transient) at every point of $M$.

Define the \emph{potential measure} $V$ associated to $Z$ by
$$
    V(A) = \int_0^\infty \mu_{s}(A)\,ds
$$
so that $ V(B) \in [0,\infty]$, for each $B \in {\cal B}(G)$. For
each $m \in M$, the {\it induced potential measure} on ${\cal
B}(M)$ is defined by

$$ V^{\Phi,m}: = V \circ \Phi_{m}^{-1} $$

In the sequel, we will frequently apply $V^{\Phi,m}$ to open balls
of the form $B_{r}(m) = \{p \in M \::\: d(p, m) < r\}$ for some
$r> 0$. When $m \in M$ is fixed, we will write $V^{\Phi}:=
V^{\Phi,m}$.

The following transience-recurrence dichotomy gives a
characterisation in terms of the behaviour of potential measures.

\begin{theorem} \label{dich}
If $Z = (Z(t), t \geq 0)$ is a \LP\ in a group $G$ and
$(G.M,\Phi)$ is a metrically invariant transformation group then
for fixed $m \in M$,
\begin{itemize}
\item[\textup{(1)}] $Z$ is either recurrent or transient at M.

\item[\textup{(2)}] $Z$ is recurrent if and only if
$V^{\Phi}(B_{r}(m)) = \infty$ for all $r > 0$.

\item[\textup{(3)}] $Z$ is recurrent if and only if
$\int_{0}^{\infty}1_{B_{r}(m)}(Z^{\Phi}(t))\,dt = \infty$ (a.s.),
for all $r> 0$.

\item[\textup{(4)}] $Z$ is transient if and only if
$V^{\Phi}(B_{r}(m)) < \infty$, for all $r > 0$.

\item[\textup{(5)}] $Z$ is transient if and only if
$\int_{0}^{\infty}1_{B_{r}(m)}(Z^{\Phi}(t))\,dt < \infty$ (a.s.),
for all $r> 0$.
\end{itemize}
\end{theorem}

\emph{Proof.} We omit the details as this is carried out in the
same way as the analogous proof for the case $G = \Rd$ which can
be found in \cite{Sato}, pp.\ 237--42. The main difference is that
we systematically replace the Euclidean norm $|\cdot|$ with $d(m,
\cdot)$ in all arguments. So analogues of (1), (2) and (3) are
first established for $G$-valued random walks (see also
\cite{GKR}, pp.\ 19--20). Observe that for each $h
> 0$, $(Z(nh), n \in \mathbb{Z}_{+})$ is a random walk
on $G$ since for each $\nN$,
$$
Z(nh) = Z(h).Z(h)^{-1}Z(2h)\cdots Z((n-2)h)^{-1}Z((n-1)h).
Z((n-1)h)^{-1}Z(nh),$$ is the composition of $n$ i.i.d.\
$G$-valued random variables.

Another important step in the proof which we emphasise is the
generalisation of the following result, due to Kingman \cite{King}
in the case $G = \Rd$, i.e. there exists $\gamma:\R^{+}
\rightarrow \R$ with $\lim_{\epsilon \downarrow 0}\gamma(\epsilon)
= 1$, such that for all $r, \epsilon, t > 0$,


 \bean \ds
   & &  P\left(\int_{t}^{\infty}1_{B_{2r}(m)}(Z^{\Phi}(s))\, ds >
   \epsilon\right)\\
& \geq &
   \gamma(\epsilon) P\Big(d(Z^{\Phi}(s+t), m) < r,\:\:\textrm{for some}\ s >
   0\Big).\eean
$\hfill \Box$

\bigskip



We now show that, just as in the case $G = \Rd$, the recurrence or
transience of \LP es, is related to that of certain embedded
random walks. Again the proof of this theorem follows along the
same lines as that in Sato \cite{Sato} p.\ 242. This time, we give
more of the details.

\begin{theorem} \label{rwalk}
If the $G$-valued random walk $(Z(nh), n \in \mathbb{Z}_{+})$ is
recurrent at $m \in M$ for some $h > 0$, then so is the \LP\ $Z$.
Conversely if $Z$ is a \cadlag\ $G$-valued \LP~that is recurrent
at $m \in M$, then there exists $h > 0$ such that the random walk
$(Z(nh), n \in \mathbb{Z}_{+})$ is recurrent.
\end{theorem}

\emph{Proof.} Suppose the random walk is recurrent at $m$, then
$\liminf_{n \rightarrow \infty}d(Z^{\Phi}(nh),m) = 0$ (a.s.).
Since
$$
    0\leq \liminf_{t\to\infty}d(Z^{\Phi}(t),m)
     \leq \liminf_{n\to\infty}d(Z^{\Phi}(nh),m) = 0,
$$
we see that $Z$ is recurrent.

Conversely, suppose that $Z$ is recurrent at $m$. First note that
$\sup_{s \in [0,h]}d(Z^{\Phi}(s), m) < \infty$ (a.s.) for each $h
> 0$ since $Z$ is \cadlag\ and $d$ is continuous. It
follows that there exists $h > 0, r > 0$ such that $P\left(
\sup_{s \in [0,h]}d(Z^{\Phi}(s), m) < r \right)
> \frac{1}{2}$. By the argument of
\cite{Sato}, p.\ 241 ((3) $\Rightarrow$ (4) therein), we deduce
that
$$
    P(Z^{\Phi}(nh) \in B_{3r}(m)) \geq \frac{1}{2h}
        \int_{(n-1)h}^{nh}P(Z^{\Phi}(t) \in B_{2r}(m))\,dt.
$$
From this and Theorem \ref{dich} (2),  it follows that
$\sum_{n=1}^{\infty}P(Z^{\Phi}(nh) \in B_{3r}(m)) = \infty$ and so
the random walk is recurrent, as required. \hfill $\Box$

\vspace{5pt}

We complete this section by establishing two straightforward but
useful results.

\begin{prop} \label{Recc} If $M$ is compact then every $G$-valued \cadlag~\LP~is
recurrent at $m \in M$.
\end{prop}

{\it Proof.} The mapping $p \rightarrow d(m,p)$ from $M$ to $\R$
is continuous and hence its image is compact. Consequently the
mapping $t \rightarrow d(m,Z^{\Phi}(t))$ is a.s. bounded for every
\cadlag~\LP~$(Z(t), t \geq 0)$ on $G$ and so such a process cannot
be transient. Hence it is recurrent by Theorem \ref{dich}(1).
$\hfill \Box$

\vspace{5pt}

{\it Remark.} Take $G=M$ in Proposition \ref{Recc} to be a compact
Lie group equipped with a left-invariant Riemannian metric. It
follows from Lemma 5.4 of \cite{AK} that (modulo some technical
conditions on the characteristics of $(\mu_{t}, t \geq 0)$) the
invariant measures for the Markov semigroup $(T_{t}, t \geq 0)$
coincide precisely with the cone of Haar measures.

\begin{theorem} \label{comp} Let $(G,M, \Phi)$ be a metrically complete
transformation group and $Z = (Z(t), t \geq 0)$ be a \LP\ in $G$.
If $Z^{\Phi,m}$ is transient at $m$ then $V^{\Phi,m}(K) < \infty$
for every compact set $K$ in $M$. If $(G,M, \Phi)$ is also
transitive then the converse statement holds.
\end{theorem}

{\it Proof.} Fix $m \in M$. If $Z$ is transient at $m$ then
$V^{\Phi,m}(B_{r}(p)) < \infty$ for all $p \in M, r > 0$. To see
this observe that we can always find a $u > 0$ such that $B_{r}(p)
\subseteq B_{u}(m)$.
Then by Theorem \ref{dich} (4), we have $V^{\Phi,m}(B_{r}(p)) \leq
V^{\Phi,m}(B_{u}(m)) < \infty$. Now since $M$ is complete, every
compact $K$ in $M$ is totally bounded and so we can find $N \in
\N, p_{1}, \ldots , p_{N} \in M$ and $r_{1}, \ldots, r_{N} > 0$
such that $K \subseteq \bigcup_{i=1}^{N}B_{r_{i}}(p_{i})$. Hence
$$ V^{\Phi,m}(K) \leq \sum_{i=1}^{N}V^{\Phi,m}(B_{r_{i}}(p_{i})) < \infty,$$
as was required.

Now suppose that $(G,M, \Phi)$ is transitive and assume that
$V^{\Phi,m}(K) < \infty$ for every compact set $K$ in $M$. Choose
$K$ to have non-empty interior $K^{0}$ and let $p \in K^{0}$. Then
we can find $r > 0$ such that $B_{r}(p) \subseteq K^{0}$. Hence
$V^{\Phi,m}(B_{r}(p)) \leq V^{\Phi,m}(K^{0}) \leq V^{\Phi,m}(K) <
\infty$ and so by transitivity $V^{\Phi,m}(B_{r}(m)) < \infty$.
Hence by Theorem \ref{dich} the process cannot be recurrent and so
it is transient.
  $\hfill \Box$

\section{Transience and Harmonic Transience for Bi-invariant Convolution Semigroups on Non-Compact Symmetric
Pairs}

In this section we will be concerned with weakly continuous
bi-invariant  convolution semigroups of probability measures
$(\mu_{t}, t \geq 0)$ defined on non-compact Riemannian symmetric
pairs $(G,K)$. We begin with some harmonic analysis of the
associated left-invariant Feller semigroup on $G$ and its
generator.


\subsection{Spherical Representation of  semigroups and
generators}

Let $(G,K)$ be a Riemannian symmetric pair of non-compact type so
$G$ is a connected semisimple Lie group with finite centre and $K$
is a maximal compact subgroup. The Iwasawa decomposition gives a
global diffeomorphism between $G$ and a direct product $KAN$ where
$A$ and $N$ are simply connected with $A$ being abelian and $N$
nilpotent wherein each $g \in G$ is mapped onto
$k(g)\exp(H(g))n(g)$ where $k(g) \in K, n(g) \in N$ and $H(g) \in
\texttt{a}$ which is the Lie algebra of $A$. We recall that the
{\it spherical functions} on $(G,K)$ are the unique mappings $\phi
\in C(G,\C)$ for which
\begin{equation} \label{spher}
\int_{K}\phi(\sigma k \tau)dk = \phi(\sigma)\phi(\tau),
\end{equation} for all $\sigma, \tau \in G$. We refer the reader
to \cite{Hel2} for general facts about spherical functions. In
particular it is shown therein that every spherical function on
$G$ is of the form
\begin{equation} \label{noncsph}
\phi_{\lambda}(\sigma) = \int_{K}e^{(i\lambda +
\rho)(H(k\sigma))}dk,
\end{equation}
for $\sigma \in G$ where $\lambda \in \texttt{a}^{*}_{\C}$ which
is the complexification of the dual space $\texttt{a}^{*}$ of
$\texttt{a}$ and $\rho$ is half the sum of positive roots
(relative to a fixed lexicographic ordering.) Note that if
$\lambda \in \texttt{a}^{*}$ then $\phi_{\lambda}$ is positive
definite.

A Borel measure $\mu$ on $G$ is said to be $K$-bi-invariant if
$$ \mu(k_{1}Ak_{2}) = \mu(A),$$
for all $k_{1}, k_{2} \in K, A \in {\cal B}(G)$. The set ${\cal
M}(K \backslash G /K)$ of all $K$-bi-invariant Borel probability
measures on $G$ forms a {\it commutative} monoid under
convolution, i.e.
$$ \mu * \nu = \nu * \mu,$$ for all $\mu, \nu \in {\cal M}(K \backslash G
/K)$. The {\it spherical transform} of $\mu \in {\cal M}(K
\backslash G /K)$ is defined by
$$   \widehat{\mu}(\lambda)  = \int_{G} \phi_{\lambda}(\sigma)\mu(d
\sigma),$$ for all $\lambda \in \texttt{a}^{*}_{\C}$. Note that

\begin{equation}  \label{stmul}
\widehat{\mu * \nu}(\lambda) =
\widehat{\mu}(\lambda)\widehat{\nu}(\lambda), \end{equation} for
all $\mu, \nu \in {\cal M}(K \backslash G /K), \lambda \in
\texttt{a}^{*}_{\C}$.

Fix a basis $(X_{j}, 1 \leq j \leq n)$ of the Lie algebra $\g$ of $G$, where $n$ is the
dimension of $G$.
~As is shown in \cite{Hunt}, \cite{Liao} there exist functions
${x_{i} \in C_{c}^{\infty}(G), 1 \leq i \leq n}$ so that $x_{i}(e)
= 0, X_{i}x_{j}(e) = \delta_{ij}$ and $\{x_{i}, 1 \leq i \leq n\}$
are canonical co-ordinates for $G$ in a neighbourhood of the
identity in $G$. A measure $\nu$ defined on ${\cal B}(G - \{e\})$
is called a {\it L\'{e}vy measure} whenever
$\int_{G-\{e\}}\left\{\left(\sum_{i=1}^{n}x_{i}(\sigma)^{2}\right)
\wedge 1\right\} \nu(d\sigma) < \infty.$

From now on we will assume that the measures forming the
convolution semigroup $(\mu_{t}, t \geq 0)$ are $K$-bi-invariant
for each $t > 0$. It then follows that for each $t > 0$, $T_{t}$
preserves the real Hilbert space $L^{2}(K \backslash G  /K)$ of
$K$-bi-invariant square integrable functions on $G$.

In this case we have Gangolli's L\'{e}vy-Khintchine formula ( see
\cite{Gang}, \cite{LW}, \cite{App0})
\begin{equation} \label{ncGLK1}
\widehat{\mu_{t}}(\lambda) = \exp\{-t\eta_{\lambda}\},
\end{equation}
for all $\lambda \in \texttt{a}^{*}, t \geq 0$. Here
\begin{equation} \label{ncGLK2}
\eta_{\lambda} = \beta_{\lambda} + \int_{G-\{e\}}(1 -
\phi_{\lambda}(\tau))\nu(d\tau), \end{equation} where
$\beta_{\lambda} \in \C$ and $\nu$ is a $K$-bi-invariant L\'{e}vy
measure on $G$ (see \cite{Gang} and \cite{LW} for details). We
call $\eta_{\lambda}$ the {\it characteristic exponent} of the
convolution semigroup. We will always assume that $(G,K)$ is {\it
irreducible}, i.e. that the adjoint action of $K$ leaves no proper
subspace of ${\textbf p}$ invariant. In this case $\beta_{\lambda}
\geq 0$. We can and will equip $G$ with a left-invariant
Riemannian metric that is also right $K$-invariant. Let
$\Delta_{G}$ be the associated Laplace-Beltrami operator of $G$.
Then $-\beta_{\lambda}$ is an eigenvalue of $a\Delta_{G}$ where $a
\geq 0$. Specifically for each $\lambda \in \texttt{a}^{*}_{\C}$
we have
\begin{equation} \label{Lapeig}
\beta_{\lambda} = a(|\lambda|^{2} + |\rho|^{2}),
\end{equation}
where the norm is that induced on $\texttt{a}^{*}_{\C}$ by the
Killing form (see \cite{Hel2}, p.427). We define
$\widetilde{\beta_{\lambda}} = |\lambda|^{2} + |\rho|^{2}$.

\begin{prop} \label{etaest}There exists $K > 0$ such that for all $\lambda \in
\texttt{a}^{*}$,
$$ |\eta_{\lambda}| \leq K(1 + |\lambda|^{2} + |\rho|^{2}).$$
\end{prop}

{\it Proof.} Let $U$ be a co-ordinate neighborhood of $e$ in $G$
and write
$$\int_{G-\{e\}}(1 -
\phi_{\lambda}(\tau))\nu(d\tau) = \int_{U-\{e\}}(1 -
\phi_{\lambda}(\tau))\nu(d\tau) + \int_{U^{c}}(1 -
\phi_{\lambda}(\tau))\nu(d\tau).$$ Since $\phi_{\lambda}$ is
positive definite we have $\sup_{\sigma \in
G}|\phi_{\lambda}(\sigma)| \leq \phi_{\lambda}(e) = 1$ and so
$\int_{U^{c}}|1 - \phi_{\lambda}(\tau)|\nu(d\tau) \leq 2\nu(U^{c})
< \infty$. Using a Taylor series expansion as in \cite{Liao} p.13
and the fact that $X_{i}\phi_{\lambda}(e) = 0$ for each $1 \leq i
\leq n$, we see that for each $\tau \in U$ there exists $\tau^{'}
\in U$ such that
$$ 1 - \phi_{\lambda}(\tau) =
-\frac{1}{2}x^{i}(\tau)x^{j}(\tau)X_{i}X_{j}\phi_{\lambda}(\tau^{\prime}).$$
Arguing as in the proof of Lemma 1 in \cite{LW}, we apply a
Schauder estimate to the equation $\Delta_{G}\phi_{\lambda} = -
\widetilde{\beta_{\lambda}}\phi_{\lambda}$ to deduce that for each
$1 \leq i,j \leq n$ there exists $C_{ij} > 0$ so that \bean
|X_{i}X_{j}\phi_{\lambda}(\tau^{\prime})| & \leq & C_{ij}(1 +
\widetilde{\beta_{\lambda}})\sup_{\sigma \in G}|\phi_{\lambda}(\sigma)|\\
& \leq &  C_{ij}(1 + \widetilde{\beta_{\lambda}}). \eean

Hence we have \bean \int_{U-\{e\}}|1 -
\phi_{\lambda}(\tau))|\nu(d\tau) & \leq & \frac{1 +
\widetilde{\beta_{\lambda}}}{2}\sum_{i,j=1}^{n}C_{ij}\int_{U-\{e\}} |x^{i}(\tau)x^{j}(\tau)|\nu(d\tau)\\
& \leq & \frac{1 +
\widetilde{\beta_{\lambda}}}{2}\left(\sum_{i,j=1}^{n}C_{ij}^{2}\right)^{\frac{1}{2}}
\int_{U-\{e\}} \sum_{i=1}^{n}x^{i}(\tau)^{2}\nu(d\tau) < \infty
\eean
and the required result follows easily from here. $\Box$

\vspace{5pt}

If $\lambda \in \at$ and $(\mu_{t}, t \geq 0)$ is symmetric, it is
easily verified that $\eta_{\lambda}$ is real-valued and
non-negative. Indeed from (\ref{ncGLK2}) and (\ref{noncsph}) we
obtain
\begin{equation} \label{ncGLK3}
\eta_{\lambda} = \beta_{\lambda} + \int_{G-\{e\}}\int_{K}\{1 -
\cos((\lambda + \rho)(H(k\tau)))\}dk \nu(d\tau).
\end{equation}

Let $C_{c}(K \backslash G /K)$ denote the subspace of $C_{c}(G)$
which comprises $K$-bi-invariant functions. If $f \in C_{c}(K
\backslash G /K)$ its spherical transform is the mapping
$\widehat{f}: \texttt{a}^{*}_{\C} \rightarrow \C$ defined by

\begin{equation} \label{sphert}
\widehat{f}(\lambda) =
\int_{G}f(\sigma)\phi_{-\lambda}(\sigma)d\sigma.
\end{equation}

We have the key Paley-Wiener type estimate that for each $N \in
\Z$ there exists $C_{N} > 0$ such that
\begin{equation} \label{spherest}
|\widehat{f}(\lambda)| \leq C_{N}(1 + |\lambda|^{-N})e^{R
|\Im(\lambda)|}
\end{equation}
where $R > 0$ (see \cite{Hel2}, p.450).

We define the Plancherel measure $\omega$ on $\texttt{a}^{*}$ by
the prescription
$$ \omega(d\lambda) = \kappa |c(\lambda)|^{-2}d\lambda,$$
where $c:\texttt{a}^{*}_{\C} \rightarrow \C$ is Harish-Chandra's
$c$-function. We will not require the precise definition of $c$
nor the value of the constant $\kappa > 0$ however we will find a
use for the estimate

\begin{equation} \label{cestimate}
  |c(\lambda)|^{-1} \leq C_{1} + C_{2}|\lambda|^{p}
\end{equation}
for all $\lambda \in \texttt{a}^{*}$ where $C_{1}, C_{2} > 0$ and
$2p = \mbox{dim}(N)$. This result follows from Proposition 7.2 in
\cite{Hel2}, p.450 and equation (16) therein on page 451.

By Theorem 7.5 of \cite{Hel2}, p.454 we have the Fourier inversion
formula for $f \in C_{c}^{\infty}(K \backslash G /K)$

\begin{equation} \label{Fourin}
f(\sigma) =
\int_{\texttt{a}^{*}}\widehat{f}(\lambda)\phi_{\lambda}(\sigma)\omega(d\lambda),
\end{equation}

which holds for all $\sigma \in G$, and the Plancherel formula:
\begin{equation} \label{Planc}
\int_{G}|f(\sigma)|^{2}d\sigma =
\int_{\texttt{a}^{*}}|\widehat{f}(\lambda)|^{2}\omega(d\lambda).
\end{equation}

By polarisation we obtain the Parseval identity for $f,g \in
C_{c}^{\infty}(K \backslash G /K)$
\begin{equation} \label{Pars}
  \la f, g \ra = \la \widehat{f}, \widehat{g} \ra,
\end{equation}

where $\la \widehat{f}, \widehat{g} \ra : =
\int_{\texttt{a}^{*}}\widehat{f}(\lambda)
\overline{\widehat{g}(\lambda)}\omega(d\lambda)$. Although the
spherical transform $f \rightarrow \widehat{f}$ extends to a
unitary transformation from $L^{2}(K \backslash G / K)$ to a
suitable space of functions on $\texttt{a}^{*}$ we cannot assume
that its precise form extends beyond the functions of compact
support. We will have more to say about this later in a special
case that will be important for us.

The next result is analogous to the pseudo-differential operator
representations obtained in Euclidean space in \cite{App1} Theorem
3.3.3 (see also \cite{App2} for extensions to compact groups.)

\begin{theorem} \label{psde}
For each $\sigma \in G, f \in C_{c}^{\infty}(K \backslash G /K)$
\begin{enumerate}
\item  $$T_{t}f(\sigma) =
\int_{\texttt{a}^{*}}\widehat{f}(\lambda)\phi_{\lambda}(\sigma)e^{-t\eta_{\lambda}}\omega(d\lambda),$$
for each $t \geq 0$. \item $${\cal A}f(\sigma) = -
\int_{\texttt{a}^{*}}\widehat{f}(\lambda)\phi_{\lambda}(\sigma)\eta_{\lambda}\omega(d\lambda).$$
\end{enumerate}
\end{theorem}

{\it Proof.}~\begin{enumerate} \item Applying Fourier inversion
(\ref{Fourin}) in (\ref{sem}) and using Fubini's theorem we obtain
\bean T_{t}f(\sigma) & = &
\int_{\texttt{a}^{*}}\int_{G}\widehat{f}(\lambda)\phi_{\lambda}(\sigma
\tau) \mu_{t}(d\tau)\omega(d\lambda)\\
& = &
\int_{\texttt{a}^{*}}\int_{G}\widehat{f}(\lambda)\int_{K}\phi_{\lambda}(\sigma
k \tau)dk \mu_{t}(d\tau)\omega(d\lambda)\\
& = & \int_{\texttt{a}^{*}}\widehat{f}(\lambda)
\phi_{\lambda}(\sigma)
\left(\int_{G}\phi_{\lambda}(\tau) \mu_{t}(d\tau)\right)\omega(d\lambda)\\
& = &
\int_{\at}\widehat{f}(\lambda)\phi_{\lambda}(\sigma)e^{-t\eta_{\lambda}}\omega(d\lambda),
\eean where we have used the left $K$-invariance of $\mu_{t}$,
(\ref{spher}) and (\ref{ncGLK1}).

\item We have \bean \frac{1}{\kappa}{\cal A}f(\sigma) & = &
\lim_{t \rightarrow 0}\int_{\at}\left(\frac{e^{-t\eta_{\lambda}} -
1}{t}\right)\widehat{f}(\lambda)\phi_{\lambda}(\sigma)|c(\lambda)|^{-2}d\lambda\\
& = & -\lim_{t \rightarrow
0}\int_{\at}e^{-\theta_{\lambda}t\eta_{\lambda}}\eta_{\lambda}\widehat{f}(\lambda)\phi_{\lambda}(\sigma)|c(\lambda)|^{-2}d\lambda,
\eean where $0 < \theta_{\lambda} < 1$ for each $\lambda \in \at,
t > 0$.

The required result follows by Lebesgue's dominated convergence
theorem. To see this first observe that for all $\lambda \in \at,
t \geq 0$, since $\Re(\eta_{\lambda}) \geq 0$, we have $
|e^{-t\theta_{\lambda}\eta_{\lambda}}| \leq 1$ and so by
Proposition \ref{etaest} for each $t \geq 0$
$$ \int_{\at}|\eta_{\lambda}e^{-t\theta_{\lambda}\eta_{\lambda}}\widehat{f}(\lambda)\phi_{\lambda}(\sigma)||c(\lambda)|^{-2}d\lambda
\leq K \int_{\at}(1+|\lambda|^{2} +
|\rho|^{2})\widehat{f}(\lambda)\phi_{\lambda}(\sigma)||c(\lambda)|^{-2}d\lambda.$$
The integral on the right hand side is easily seen to be finite by
using (\ref{cestimate}) and taking $N$ to be sufficiently large in
(\ref{spherest}).  $\hfill \Box$

\end{enumerate}


{\bf Note}~Based on the results of Theorem \ref{psde} we may write
$$ \widehat{T_{t}f}(\lambda) = e^{-t
\eta_{\lambda}}\widehat{f}(\lambda)~~\mbox{and}~~\widehat{{\cal
A}f}(\lambda) = - \eta_{\lambda}\widehat{f}(\lambda),$$ for all $t
\geq 0, \lambda \in \at$. 


We will need to extend the precise form of Parseval's formula to
include the range of $T_{t}$ acting on $C_{c}(K \backslash G /K)$.

\begin{theorem}  \label{Parmore}
For all $f,g \in C_{c}(K \backslash G /K), t \geq 0$,
\begin{equation} \label{Parmore1} \la T_{t}f, g \ra = \la
\widehat{T_{t}f}, \widehat{g} \ra. \end{equation}
\end{theorem}

{\it Proof.} Using the result of Theorem \ref{psde}(1), Fubini's
theorem, the fact that for all $\sigma \in G, \lambda \in \at,
\overline{\phi_{\lambda}(\sigma)} = \phi_{-\lambda}(\sigma)$ and
(\ref{sphert}), we find that for all $f,g \in C_{c}(K \backslash G
/K), t \geq 0$, \bean \la T_{t}f, g \ra & = &
\int_{G}\left(\int_{\texttt{a}^{*}}\widehat{f}(\lambda)\phi_{\lambda}(\sigma)e^{-t\eta_{\lambda}}\omega(d\lambda)\right)
g(\sigma)d\sigma \\
& = &
\int_{\texttt{a}^{*}}e^{-t\eta_{\lambda}}\widehat{f}(\lambda)
\left(\overline{\int_{G}g(\sigma)\phi_{-\lambda}(\sigma)d\sigma}\right)\omega(d\lambda)\\
& = &
\int_{\texttt{a}^{*}}e^{-t\eta_{\lambda}}\widehat{f}(\lambda)\overline{\widehat{g}(\lambda)}\omega(d\lambda)\\
& = & \la \widehat{T_{t}f}, \widehat{g} \ra
~~~~~~~~~~~~~~~~~~~~\hfill \Box \eean

Now assume that the convolution semigroup $(\mu_{t}, t \geq 0)$ is
symmetric, i.e. $\mu_{t} = \widetilde{\mu_{t}}$ for all $t \geq
0$. The space $C_{c}^{\infty}(K \backslash G /K)$ is a dense
linear subspace of the Hilbert space $L^{2}(K \backslash G /K)$ of
$K$-bi-invariant square integrable functions on $G$. We consider
the restriction therein of the Dirichlet form ${\cal E}$ where
${\cal E}(f): = ||(-{\cal A})^{\frac{1}{2}}f||^{2}$ on ${\cal D}:=
\mbox{Dom}({\cal A}^{\frac{1}{2}})$.

\begin{cor} \label{invDir}
For each $f, g \in C_{c}^{\infty}(K \backslash G /K)$ we have
$$ {\cal E}(f,g) =
\int_{\texttt{a}^{*}}\widehat{f}(\lambda)\eta_{\lambda}\overline{\widehat{g}(\lambda)}\omega(d\lambda).$$
\end{cor}

{\it Proof.} This follows from  differentiating both sides of
(\ref{Parmore1}) using Theorem \ref{psde}(1). $\hfill \Box$

\subsection{Harmonic Transience}

Let $(\mu_{t}, t \geq 0)$ be an arbitrary bi-invariant convolution
semigroup in $G$. We define the {\it potential operator} $N$ on
$C_{0}(G/K)$ by the prescription

$$ Nf: = \lim_{t \rightarrow \infty}\int_{0}^{t} S_{u}^{\Phi}f du$$

for all $f \in \mbox{Dom}(N):= \left\{f \in C_{0}(G/K); \lim_{t
\rightarrow \infty}\int_{0}^{t} S_{u}^{\Phi}f du~\mbox{exists in}~
C_{0}(G/K)\right\}$. We say that $(\mu_{t}^{\pi}, t \geq 0)$ is
{\it integrable} if $C_{c}(G/K) \subseteq \mbox{Dom}(N)$. By the
same arguments as are presented in Lemmata 12.2 and 13.19 and
Proposition 13.21 of \cite{BF}, it follows that integrability
implies transience. Berg and Farault \cite{BeFa} have shown that
if $G/K$ is irreducible then the projection to $G/K$ of every
bi-invariant convolution semigroup on $G$ is integrable and hence
all of these semigroups of measures are transient. It then follows
that $({\cal E}, D)$ is a transient Dirichlet space. In Theorem
7.3 in \cite{Hey2} this result is extended to the case where the
symmetric space is no longer required to be irreducible.
Furthermore it is shown that the the associated Feller semigroup
is {\it mean ergodic}, i.e. $\lim_{T \rightarrow
\infty}\frac{1}{T}\int_{0}^{T}S_{u}^{\phi}fdu$ exists for all $f
\in C_{0}(G/K)$.

We say that $(\mu_{t}, t \geq 0)$ is {\it harmonically transient}
if the mapping $$\lambda \rightarrow
\ds\frac{1}{\Re(\eta_{\lambda})}$$ from $\at$ to $\R$ is locally
integrable with respect to the Plancherel measure $\omega$.

Note that the result of (i) in Theorem \ref{trans} below is
well-known to be a necessary and sufficent condition for a
transient Dirichlet space (see \cite{Deny}) however we include a
short proof to make the paper more self-contained.

\begin{theorem} \label{trans} If $(G,K)$ is an irreducible Riemannian symmetric pair and
$(\mu_{t}, t \geq 0)$ is a symmetric $K$-bi-invariant convolution
semigroup then
\begin{enumerate}
\item[(i)] $\int_{0}^{\infty}\la T_{t}f, f \ra dt < \infty$ for
all $f \in C_{c,+}(K \backslash G / K)$. \item[(ii)] $(\mu_{t}, t
\geq 0)$ is harmonically transient.
\end{enumerate} \end{theorem}

{\it Proof.}~\begin{enumerate} \item[(i)] Let $f \in C_{c,+}(K
\backslash G /K)$ and let $C = \mbox{supp}(f)$ then
$$\int_{0}^{\infty}\la T_{t}f, f \ra dt =
\int_{C}f(\sigma)\int_{G}1_{\sigma^{-1}C}(\tau)f(\sigma\tau)V(d\tau)d\sigma.$$
Define $C^{-1}C := \{\sigma^{-1}\tau; \sigma,\tau \in C\}$. Since
$C^{-1}C$ is the image of the continuous mapping from $C \times C$
to $G$ which takes $(\sigma,\tau)$ to $\sigma^{-1}\tau$ it is
compact. It follows that $\sup_{\sigma \in C}V(\sigma^{-1}C) \leq
V(C^{-1}C) < \infty$ as $(\mu_{t}, t \geq 0)$ is transient, and so
the mapping $\sigma \rightarrow
\int_{G}1_{\sigma^{-1}C}(\tau)f(\sigma\tau)V(d\tau)$ is bounded.
The required result follows.

\item[(ii)] This is proved by a similar method to the Euclidean
space case (see Example 1.5.2 in \cite{FOT} p.42-3. ) We include
it for the convenience of the reader, noting that it was for this
specific purpose that we established Theorem \ref{Parmore}. 
We choose $f \in C_{c,+}(K \backslash G /K)$ such that
$\widehat{f}(\lambda) \geq 1$ for all $\lambda$ in some compact
neighbourhood $A$ of $\at$. Then using Theorems \ref{psde} (1) and
\ref{Parmore} we get \bean
\int_{A}\frac{\omega(d\lambda)}{\eta_{\lambda}} & \leq &
\int_{A}\frac{|\widehat{f}(\lambda)|^{2}}{\eta_{\lambda}}\omega(d\lambda)
 \leq  \int_{\at}\frac{|\widehat{f}(\lambda)|^{2}}{\eta_{\lambda}}\omega(d\lambda)\\
& = & \int_{0}^{\infty}\la \widehat{T_{t}f}, \widehat{f} \ra dt =
 \int_{0}^{\infty}\la T_{t}f, f \ra dt  < \infty.\eean   $\hfill \Box$

\end{enumerate}

\begin{cor} If $(G,K)$ is an irreducible Riemannian symmetric pair and
$(\mu_{t}, t \geq 0)$ is a $K$-bi-invariant convolution semigroup
then it is harmonically transient.
\end{cor}

{\it Proof.} Let $(\mu_{t}, t \geq 0)$ be an arbitrary
bi-invariant convolution semigroup with characteristic exponent
$\eta_{\lambda}$. We associate to it a symmetric bi-invariant
convolution semigroup $(\mu_{t}^{(R)}, t \geq 0)$ by the
prescription $\mu_{t}^{(R)} = \mu_{t} * \widetilde{\mu_{t}}$. It
follows easily from (\ref{stmul}) and (\ref{ncGLK1}) that it has
characteristic exponent $2 \Re(\eta_{\lambda})$ and using this
fact we can immediately deduce that $(\mu_{t}, t \geq 0)$ is
harmonically transitive by using Theorem \ref{trans}. $\hfill
\Box$

\vspace{5pt}

{\bf Acknowledgements.}~I am grateful to Herbert Heyer, Ming Liao,
Ren\'{e} Schilling and Thomas Simon for useful discussions.


\begin{thebibliography}{99}

\bibitem{App0} D.Applebaum, Compound Poisson processes and L\'{e}vy
processes in groups and symmetric spaces, {\it J.Theor. Prob.}
{\bf 13}, 383-425 (2000)

\bibitem{App1} D.Applebaum, {\it L\'{e}vy Processes and
Stochastic Calculus} (second edition), Cambridge University Press
(2009)

\bibitem{App2} D.Applebaum,  Infinitely divisible central probability measures on compact
Lie groups - regularity, semigroups and transition kernels, to
appear in {\it Annals of Prob.} (2011)


\bibitem{AK} D.Applebaum, H.Kunita, Invariant measures for
L\'{e}vy flows of diffeomorphisms, {\it Proc. Roy. Soc. Edinburgh}
{\bf 130A}, 925-46 (2000)

\bibitem{Bal} P.Baldi, Caract\'{e}risation des groupes de Lie connexes
r\'{e}current, \emph{Ann.\ Inst.\ H.\ Poincar\'{e} (Probab.\
Stat.)} {\bf 17}, 281-308 (1981).

\bibitem{BLP} P.Baldi, N.Lohou\'{e}, J.Peyri\`{e}re, Sur la
classification des groupes r\'{e}currents, {\it C.R.Acad.Sc.Paris,
S\'{e}rie A} {\bf t.285}, 1103-4 (1977)

\bibitem{Berg} C.Berg, Dirichlet forms on symmetric spaces, {\it
Ann. Inst. Fourier, Grenoble} {\bf 23}, 135-56 (1973)

\bibitem{BeFa} C.Berg, J.Faraut, Semi-groupes de Feller invariants
sur les espaces homog\`{e}nes non moyennables, {\it Math. Z.} {\bf
136} 279-90 (1974)

\bibitem{BF} C. Berg, G. Forst, {\it Potential Theory on Locally
Compact Abelian Groups}, Springer-Verlag (1975).

\bibitem{BH} W. R. Bloom, H. Heyer, {\it Harmonic Analysis of
Probability Measures on Hypergroups}, de Gruyter (1995).


\bibitem{Born} E.Born, An explicit L\'{e}vy-Hin\v{c}in formula for convolution
semigroups on locally compact groups, {\it J.Theor. Prob.} {\bf 2}
325-42 (1989)

\bibitem{Chav} I.\ Chavel, \emph{Riemannian Geometry: A Modern
Introduction}, Cambridge University Press, Cambridge (1993).

\bibitem{Deny} J. Deny, M\'{e}thodes Hilbertiennes et th\'{e}orie
du potentiel. In {\it Potential Theory}, ed.  M. Brelot, Centro
Internazionale Matematico Estivo, Edizioni Cremonese, pp.~123--201
(1970), reprinted in {\it Potential Theory - CIME Summer Schools}
{\bf 49}, Springer-Verlag Berlin Heidelberg (2010)

\bibitem{EK} S. N. Ethier, T. G. Kurtz, {\it Markov Processes,
Characterisation and Convergence}, Wiley (1986).


\bibitem{FOT} M. Fukushima, Y. Oshima, M. Takeda, {\it Dirichlet
Forms and Symmetric Markov Processes}, de Gruyter (1994).


\bibitem{Gang} R.Gangolli, Isotropic infinitely divisible measures on
symmetric spaces, {\it Acta Math}, {\bf 111}, 213-46 (1964)

\bibitem{Gr} A.\ Grigor'yan, Analytic and geometric background of recurrence and
non-explosion of the Brownian motion on Riemannian manifolds,
\emph{Bull.\ Am.\ Math.\ Soc.} {\bf 36}, 135-249 (1999).

\bibitem{GKR} Y.\ Guivarc'h, M.\ Keane, B.\ Roynette, \emph{Marches
Al\'{e}atoires sur les Groupes de Lie}, Lect. Notes Math. {\bf
624}, Springer, Berlin (1977).

\bibitem{Hel1} S.Helgason, {\it Differential Geometry, Lie Groups and
Symmetric Spaces}, Academic Press (1978) \vspace{5pt}

\bibitem{Hel2} S.Helgason, {\it Groups and Geometric Analysis}, Academic
Press (1984)

\bibitem{Hey} H.Heyer, {\it Structural Aspects in the Theory of
Probability}, World Scientific (2005)

\bibitem{Hey1} H.Heyer, Convolution semigroups of probability measures on Gelfand pairs,
{\it Expo. Math.} {\bf 1}, 3-45 (1983)

\bibitem{Hey2} H.Heyer, Transient Feller semigroups on certain
Gelfand pairs, {\it Bull. Inst. Mat. Acad. Sin.} {\bf 11}, 227-56
(1983)

\bibitem{HT} H.Heyer, G.Turnwald, On local integrability and
transience of convolution semigroups of measures, {\it Acta Appl.
Math.} {\bf 18}, 283-96 (1990)

\bibitem{Hunt} G.A.Hunt, Semigroups of measures on Lie groups, {\it Trans. Amer.
Math. Soc.} {\bf 81}, 264-93 (1956)

\bibitem{It} M.It\^{o}, Transient Markov convolution semi-groups
and the associated negative definite functions, {\it Nagoya Math.
J.} {\bf 92}, 153-61 (1983); Remarks on that paper, {\it Nagoya
Math. J.} {\bf 102}, 181-84 (1986)

\bibitem{King} J.F.C.Kingman, Recurrence properties of processes
with stationary independent increments, {\it J. Austral. Math.
Soc.} {\bf 4}, 223-8 (1964)

\bibitem{Liao} M.Liao, {\it L\'{e}vy Processes in Lie
Groups}, Cambridge University Press (2004)

\bibitem{LW} M.Liao, L.Wang, L\'{e}vy-Khinchine formula and
existence of densities for convolution semigroups on symmetric
spaces, {\it Potential Analysis} {\bf 27}, 133-50 (2007)

\bibitem{PS} S.C.\ Port, C.J.\ Stone, Infinitely divisible processes
and their potential theory II, \emph{Ann.\ Inst.\ Fourier} {\bf
21}.4, 179-265 (1971)

\bibitem{Sato} K-I.Sato, {\it L\'{e}vy Processes and Infinite
Divisibility}, Cambridge University Press (1999)

\bibitem{Wo} J.A.Wolf, {\it Harmonic Analysis on Commutative
Spaces} Amer. Math. Soc. (2007)



\end{thebibliography}
\end{document}